\documentclass[11pt]{amsart}
\usepackage{amsfonts,amsmath,amssymb}
\usepackage[all]{xy}
\usepackage{hyperref}
\begin{document}

\newcommand{\ba}{{\bf a}}
\newcommand{\bb}{{\bf b}}
\newcommand{\fb}{{\sf b}}
\newcommand{\ab}{{\rm ab}}
\newcommand{\bpi}{{\boldsymbol{\pi}}}
\newcommand{\be}{{\bf e}}
\newcommand{\oh}{{\mathfrak o}}
\newcommand{\m}{{\mathfrak m}}
\newcommand{\jnf}{{\rm inf}}
\newcommand{\A}{{\mathbb A}}
\newcommand{\C}{{\mathbb C}}
\newcommand{\F}{{\mathbb F}}
\newcommand{\G}{{\mathbb G}}
\newcommand{\M}{{\mathbb M}}
\newcommand{\N}{{\mathbb N}}
\newcommand{\bP}{{\mathbb P}}
\newcommand{\R}{{\mathbb R}}
\newcommand{\bS}{{\mathbb S}}
\newcommand{\Gal}{{\rm Gal}}
\newcommand{\Q}{{\mathbb Q}}
\newcommand{\bQ}{{\overline{\mathbb Q}}}
\newcommand{\T}{{\mathbb T}}
\newcommand{\Sp}{{\rm Sp \;}}
\newcommand{\W}{{\sf W}}
\newcommand{\sK}{{\sf K}}
\newcommand{\Hom}{{\rm Hom}}
\newcommand{\Z}{{\mathbb Z}} 
\newcommand{\Be}{{\boldsymbol{\epsilon}}}

\newcommand{\sslash}{{/\!/}}

\newcommand{\BHV}{{\rm BHV}}
\newcommand{\Spec}{{\rm Spec \; }}

\parindent=0pt
\parskip=6pt

\newcommand{\ie}{\textit{ie}\,}
\newcommand{\eg}{\textit{eg}\,}
\newcommand{\se}{{\sf e}}
\newcommand{\cf}{{\textit{cf}\,}}

\newcommand{\essl}{{\mathfrak {sl}}}
\newcommand{\Rep}{{\rm Rep}}
\newcommand{\ox}{{\rm or}}
\newcommand{\bM}{{\overline{\mathbb M}}}
\newcommand{\Mod}{{\rm Mod}}
\newcommand{\Cr}{{\rm cr}}

% version of 1 June 025
\author[J Morava]{J Morava}

\address{Department of Mathematics, The Johns Hopkins University,
Baltimore, Maryland} 

\email{jmorava1@jhu.edu}

\title{On a topological complex orientation for $\T$-equivariant K-theory} 

\begin{abstract}{The principal result of this note is the existence of a complex topological orientation for Atiyah-Segal $\T$-equivariant K-theory which indexes the space $\C P^n$ of lines in $\C^{n+1}$ by the Fourier expansion $1 + q + \dots + q^n$}\end{abstract}
  
\maketitle

Acknowledgement: This work began in conversations with Markus Szymik at the Mittag-Leffler Institute, February 2017. It is a pleasure to thank them for support, and him, together with Andrew Baker, Sanath Devalapurkar, and Clarence Wilkerson, for conversations about this and related topics. \bigskip

{\bf \S 1 Projective varieties and complex cobordism} \bigskip

{\bf 1.0} The characteristic polynomial \cite{6}
\[
\chi_{Y,Z}(M) = \sum [H^{i,j}(M)] Y^i Z^j \in K^0 ({\rm pt})[Y,Z]
\]
($\chi_{Y,Z} = \chi_{Z,Y}$ by Serre duality) defined by the Hodge-Dolbeault cohomology of a smooth complex projective variety has several specializations of topological interest \begin{footnote}{\eg $\chi_{iY,iY}(M) = \chi(M) Y^{{\rm dim} M}$ recovers the graded Euler characteristic}\end{footnote}. It is natural and multiplicative  
\[
\chi_{Y,Z}(M_0 \times M_1) = \chi_{Y,Z}(M_0) \cdot \chi_{Y,Z}(M_1) 
\]
with respect to Cartesian product, defining a ring homomorphism
\[
\chi_{Y,Z} : \{{\rm Proj {\;} Var}\} \to  K^0({\rm pt})[Y,Z] 
\]
which, for example, sends projective space to 
\[
\chi_{Y,Z}(\C P^n)  = 1 + YZ + \cdots + (YZ)^n {\;}. \bigskip
\]
{\bf 1.1} Following Serre \cite{18}(Ch IV \S 5 remarks) K\"ahler's identities on the Hodge cohomology \cite{12}(appendix 2) of a smooth projective variety define a natural action of the Lie algebra $\essl_2(\C)$; see \cite{16}(appendix 2) for a symplectic analog. The Grothendieck ring $\Rep(\essl_2)$  of such things is a free $\lambda$-algebra on one generator $\sigma$.

On another hand, Quillen's construction of complex cobordism  assigns to a compact (stably almost complex $\dots$) $2n$-dimensional manifold $X$ with an action of the circle $\T$, an element
\[ 
[X \sslash \T]  = [X] + \Sigma_{i \geq 1}  [X_i] c^i  \in MU^{-2n} B\T, {\;} [X_i] \in MU_{2(n+i)}
\]
of the complex cobordism ring of the classifying space $B\T \simeq \C P^\infty$ (where $c$ is Quillen's Euler-Chern class \cite{13} in $MU^2 B\T$); there are details below in an appendix.

{\bf 1.2} In \S 2 we define a formal group law
\[
X +_\chi Y  = F_\chi(X,Y) := \frac{X + Y + (1+ q)XY}{1 + qXY} \in K_\T [[X,Y]]
\]
over $\Z[q] \subset \Z[q,q^{-1}] \cong K^0_\T({\rm pt})$, regarded as the Atiyah-Segal equivariant $\Z_2$ - graded K - theory of a point with circle action \cite{2}. The classifying homomorphism \cite{3}
\[
MU^*(B\T) \ni c \to \exp_\chi (c) \in K_\T(B\T)
\]
of graded rings defines a complex topological orientation for $K_\T$, sending
\[
MU^{-2n} \ni [\C P^n] \to [n]_q = 1 + q + \cdots + q^n \in K_\T {\;}.
\]

{\bf Hypothesis :} {\it There is a commutative diagram 
\[
\xymatrix{
\{{\rm Proj {\;} Var}\} \ar[dr] \ar[r] & \Rep(\essl_2) \ni \sigma \ar@{.>}[d]^? \ar@{.>}[dr] \ar[r] & \chi_{Y,Z}(\sigma) \in \Z[Y,Z] \ar@{.>}[d]^{YZ \mapsto q} \\
{ } & MU^* \ar[r]^-\chi & K_\T = \Z[q,q^{-1}]}
\]
with a morphism of $\lambda$-rings along the diagonal.}

{\bf motivation :} This suggestion originated with issues involved in extending previous work \cite{10} on Swan-Tate $\T$-equivariant cohomology $t^*_\T H$ to Swan-Tate $\T$-equivariant K-theory $t^*_\T K$. It seems to say something about the global homotopy theory of projective varieties and perhaps about metaplectic representations \cite{11} of the Virasoro group. 

The principal {\bf result} of this note is the existence of a complex orientation for Atiyah-Segal $\T$-equivariant K-theory which indexes the space $\C P^n$ of lines in 
$\C^{n+1}$ by the Fourier expansion $1 + q + \dots + q^n$. \bigskip

{\bf 1.3} Let us be more precise. The quasi-character $q : \T \subset \C^\times$ identifies the $K$ - theory of vector bundles 
\[
K^0([{\rm pt} \sslash \T]) = K_\T (\rm pt) \cong \Z[q,q^{-1}]
\]
over a point with $\T$-action, with the representation ring $\Rep_\C(\T)$ \begin{footnote}{with $\T \ni u \mapsto ((u,z) \to u^n z)$ for the representation $[\C(n)]$}\end{footnote}. This is a (special, classical) $\lambda$-ring, with Adams operations $\psi^k(q) = q^k$.
 It is natural to identify the localization defined by 
\[
t_\T^* MU \to t_\T K \cong (1-q)^{-1}K_\T  \in \Z[q,q^{-1},(1-q)^{-1}] - \Mod
\]
with $\T$-equivariant Swan-Tate K-theory \cite{14}, see \S 3.2 below. \bigskip 

{\bf \S 2 construction of} $F_\chi$ \bigskip

{\bf 2.1 lemma :} {\it If

\[
[Q] \; = \; \left|\begin{array}{cc} - q & 1 \\
                        -1  & 1 \end{array}\right| ,  \;                                                                                   
[Q^{-1}] \; = \; \left|\begin{array}{cc} 1 & - 1 \\
                           1 & - q \end{array}\right| 
\]
then }
\[
[Q] \circ [Q^{-1}] = [Q^{-1}] \circ [Q] \; = \; (1 - q)
\left|\begin{array}{cc} 1 & 0 \\
                        0 & 1 \end{array}\right| \; .
\]

Note that $\det \; [Q] = \det \; [Q^{-1}] = (1 - q) := |Q|$. $\Box$\bigskip

Recall that a $2 \times 2$ matrix 
\[
[A] := \left|\begin{array}{cc} a & b \\
                        c & d \end{array}\right| \; .
\]
defines
\[
[A](T) := \frac{aT + b}{cT + d} \;,
\] 
\eg $[Q](T) = 1 + (1-q)T + \dots$. \bigskip

{\bf 2.2 definition} 
\[
\log_\chi(T) := |Q|^{-1} \log \: [Q](T) = T + \dots \in \Q[q][[T]]
\]
and 
\[
\exp_\chi(T) :=  [Q^{-1}](\exp \: |Q| T)  = T + \dots \in \Q[q][[T]]
\]

satisfy $(\exp_\chi \circ \log_\chi)(T) = T$. \bigskip

{\bf 2.3 proposition:} 
\[
X +_\chi Y  = F_\chi(X,Y) := \frac{X + Y + (1+ q)XY}{1 + qXY} \in \Z[q][[X,Y]]
\]
{\it defines a formal group.}

{\bf proof :} Things being torsion-free, we can work over $\Q$ and use logarithms to show that 
\[
X +_\chi Y  =  \exp_\chi (\log_\chi(X) + \log_\chi(Y)) \;.
\]
A straightforward expansion of the righthand expression yields
\[
[Q]^{-1}(\exp [ (|(Q| \cdot |Q|^{-1}) [\log [Q](X) + \log [Q](Y)]] = 
\]
\[
= [Q]^{-1}([Q](X) \cdot [Q](Y)) = \frac{[Q](X) \cdot [Q](Y) - 1}{[Q](X) \cdot [Q](Y) - q} =
\]

\[
= \frac{(1-qX)(1-qY) - (1-X)(1-Y)}{(1-qX)(1-qY) - q(1-X)(1-Y)} =
\]

\[
= \frac{(1 - q(X+Y) + q^2 XY) - (1 - (X+Y) -XY)}{(1 - q(X+Y) + q^2 XY) - (q - q(X+Y) + qXY)} = 
\]

\[
= \frac{ - (q-1)(X+Y) + (q^2-1) XY}{1-q  - q(q-1)XY} =
\]

\[
= \frac{(1-q)(X+Y + (1+q)XY)}{(1-q)(1 + qXY)} 
\]
$\& \dots$ as was to be shown. \bigskip

Note that $q=0$ specializes $F_\chi$ to the multiplicative formal group law
\[
\G_m(X,Y) = X + Y + XY
\]
of arithmetic.\bigskip

$\bullet$ {\bf 2.4 exercise 1} 
\[
\log_\chi (T) = \sum_{k \geq 1} [k]_q \frac{T^k}{k}
\]

where the `$q$' - numbers 
\[
[k]_q = 1 + q + \dots + q^{k-1}  = \frac{1 - q^k}{1 - q} \in (1 + q\Z[[q]])^\times
\]
are `$q$' - adic units. \bigskip

{\bf solution:} 

Unwinding the definition, we have
\[
 |Q|^{-1} \log \: [Q](T)  = \frac{1}{1-q}[\log(1 - qT) - \log(1-T)] =  
 \]
 
 \[
 = \frac{-1}{1-q} \sum_{k \geq 1} [\frac{(qT)^k}{k} - \frac{T^k}{k}] =  
 \sum_{k \geq 1} \frac{q^k - 1}{q - 1}  \frac{T^k}{k} =  \sum_{k \geq 1} [k]_q \frac{T^k}{k} .
 \] \bigskip

 {\;} $\bullet$ 2
\[ 
q^{1/2}F_\chi (q^{-1/2}X,q^{-1/2}Y) = \frac{X + Y + (q^{-1/2} + q^{1/2})XY}{1 + XY}
\]
\cite{4}(\S 5.4.2)\bigskip

 {\;} $\bullet$ 3
\[
\exp_{\G_{\hat{m}}}( t \log_\chi(T)) = 1 - [Q](T)^{-|Q|t}{\;},
\] 
\cite{13, 15}

 {\;} $\bullet$ 4
\[
[n]_\chi(T) = \exp_\chi (n \log_\chi (T)) = [Q^{-1}]([\Q](T)^n) \in \Z[q][[T]]
\] 

 {\;} $\bullet$ 5 The translation-invariant differential 
 \[
 d\log_\chi(T) = q^{-1} [(T - q^{-1})^{-1} - ((T - 1)^{-1} ] \cdot dT = 
 \]
\[
\frac{dT}{1 - (1+ q) T + q T^2}
\]
can be imagined as a dipole merging as $q \to 1$.

{\;} $\bullet$ 6 : Reducing modulo $p$, we have
\[
[p]_\chi(T) \equiv \frac{[p]_q}{1 + [p-1]_q \cdot q T^p} {\;} T^p
\]
\bigskip
       
{\bf 2.5} It follows that $F_\chi \otimes \Z_p$ has height one away from torsion points on circle, so Landweber exactness 
\[
X \mapsto MU^*(X) \to MU^*(X) \otimes_\chi K_\T \cong K_\T(X)
\]       
defines a cohomology theory by an argument going back to Conner and Floyd. We thereby have a $\T$-equivariant lift \cite{8}(lemma 2.3)
\[
F_\chi : MU^* B\T \ni c \mapsto \exp_\chi (T) \in K_\T [[T]]
\]
of the genus $\chi$ to an equivariant orientation for circle actions on complex-oriented manifolds. 

{\bf \S 3 some extensions} 
             
{\bf definition} The localization \cite{14}(\S 3)
\[
\Z[q] \subset \Z[q]^\Cr  := \Z[q][[k]_q^{-1} \:|\: k \geq 1] \subset \Z[[q]]
\]
defines the ring of cromulent functions of $q$. For example, we write $[k]_q!^{-1}$ for the cromulent function $\prod^n_1 [i]_q^{-1}$. Note that $q \to 0$ maps  $\Z[q]^\Cr \to \Z$, while $q \to 1$ sends  $\Z[q]^\Cr \to \Q$ : Spec $\Z[q]^\Cr$ is a span 
\[
\Spec \Z \coprod \Spec \Q  \to \Spec \Z[q]^\Cr {\;},
\]
something like the spectrum of the integers, with a cone over the generic point. 

{\bf 3.1} The $\lambda$-ring structure on the representation ring $\Rep_\C(\T)$ sends a line $L$ to $1 + tL$ and is multiplicative, which suggests that 
\[
\lambda_{-t} (1 - q)^{-1}  = \lambda_{-t }(\Sigma_{k \geq 0} q^k) = \prod_{n \geq 0}\lambda_{-t}(q^n) = (t;q)_\infty
\]
(using $q$-Pochhammer notation $(t)_n := (t;q)_n = \prod_{n-1 \geq k \geq 0} (1 - tq^k)$). Per Wikipedia we have that 
\[
(t;q)_\infty = \sum_{k \geq 0} {\;} (- 1)^k [k]!_q^{-1} q^{\binom{k}{2}}{\;} (1-q)^{-k} t^k
\]
which confirms the identification
\[
\lambda^k (1- q)^{-1} = [k]!_q^{-1} q^{k(k+1)/2}(1-q)^{-k} = (-1)^k [-k]_q!^{-1} (1-q)^{-k}
\]
(the latter makes sense, via $\Gamma_{[q]} \dots$). This illustrates examples 3, 4 of MacDonald \cite{9}(Ch I \S 2.15 p 26) where it is shown that 
\[
\psi^k(1 - q)^{-1} = (1 - q^k)^{-1} 
\]
in slightly different notation. Following Kedlaya \cite{7}(example 4.2.4), \cite{14}(\S 3) we have a generalized $\lambda$-ring structure
\[
\lambda_t : (1-q)^{-1} K_\T = t_\T K \to \W(t^\vee_\T K )
\]
taking values in Swan-Tate Atiyah-Segal $\T$-equivariant K-theory
\[
t^\vee_\T K := t_\T K \otimes_{\Z_[q]} \Z[q]^\Cr
\]
with cromulent coefficients.\bigskip

{\bf 3.2} The expression variously written 
\[
\lambda_{-1}(1 - q)^{-1} =  \prod_{k \geq 1} (1 - q^k) = \phi(q) = \Sigma_\Z (-1)^n q^{n(3n-1)} = (q;q)_\infty 
\]
has a topological interpretation as a Thom class (Atiyah 67 p 100) for the inverse of the virtual vector bundle $1 - q = [ \C(0)] -  [\C(1)] \in K_\T({\rm pt})$. 

We may regard $[\phi^{\pm 1}] t_\T K$ as a formally 24-periodic height one cohomology functor supported on the zero-locus
\[
 \Delta \sim E_4^3 - E_6^2 = 0
\]
of nodal (topologically \cite{5}(Prop 5.7) $\C_+/\{0,1\}$) curves in the compactified elliptic moduli stack (where $j = \infty$, but $E_4$ is unrestricted \cite{19}). 

[IIUC this instantiates the motivic splitting 
\[
{\mathbb P}^1(\C) \sim \G_m(\R)_+ \wedge \C_+ / \{0,1\} {\;}.
\]

{\bf 3.3 exercise} 
\[ \lim_{t \to 1} q \lambda_{-t}
 \left|\begin{array}{cc} 0 & 24 \\
                                  -1  & 1 \end{array}\right| (q)
\sim \Delta(q)
\]
is the modular discriminant \cite{1}.

{\bf appendix :} {\it Quillen rules OK?} 

$\; \bullet$ A (compact closed cx-oriented) $2n$ - manifold $X$ with $\T$-action defines $[X] \in MU^{-2n} B\T$, while a ($\dots 2m$) - manifold $Y$ with a {\bf free} $\T$ - action defines $[Y] \in MU_{2m} B\T$. 

The product $X \times Y$, with diagonal $\T$ - action, defines a module structure
\[
\cap : [X], [Y] \to [X \times Y] 
\]
leading to an $MU$ -module isomorphism
\[
MU^* B\T \to {\rm Hom}^*_{MU}(MU_*B\T,MU_*) {\;}.
\]
and a commutative diagram
\[
\xymatrix{
MU^* B\T \otimes_{MU} MU_*B\T  \ar[d] \ar[r] & MU_*B\T \ar[d] \\
H^*(B\T,MU) \otimes_{MU} H_*(B\T,MU) \ar[r] & H_*(B\T,MU)} {\;.}
\]
Here $H^*B\T \cong \Z[c], {\;} H_*B\T \cong \Z[\gamma^*b] \subset \Q[b]$ is a divided power algebra, $c^i \cap \gamma^k b = \gamma^{k-i} b$ (\ie $ c = \partial_b$), and $MU_*B\T_+ \cong MU_*[b_k \dots]$ with
\[
b(T) = 1 + \Sigma_{k \geq 1} b_k T^k = \exp (b \log_{MU}(T)) = \Sigma_n {\;} \Pi_k {\;}  (\frac{\C P^{n-1}}{n})^k \cdot \gamma^I b \cdot T^n
\]
(where $\gamma^I b = \Pi_1^n (\gamma^k b)^{r_k}, {\;} | I | =  \Sigma i r_i $, \cite{10, 12}). Thus 
\[
c \cap b(T) = \log_{MU}(T) \cdot b(T) \Rightarrow
\]
\[
c \cap b_k = \Sigma {\;} j^{-1} \C P^{j-1} \cdot b_{k-j}  = b_{k-1} + \cdots + k^{-1} \C P^{k-1} {\;}.
\]

{\bf example:} If $[\C P^1(\omega)] \in MU_2 B\T_+$ represents the symplectic cobordism class defined by the Fubini metric, while $[\C P^1(0)]$ represents the underlying space, then \cite{15}
\[
b_1 = [\C P^1(\omega)] - [\C P^1(0)] = c^{-1} \in t^{-2}_\T MU \subset t_\T^*MU \cong MU^*((c)) {\;}.
\]

\bigskip

\bibliographystyle{amsplain}

\end{document}